\documentclass[twoside,12pt]{amsart}

\usepackage{amsmath}
\usepackage{amssymb}
\usepackage{amscd}
\usepackage{amsthm}
\usepackage{Kropina}

\begin{document}

\rightline{The translation from Russian}
\rightline{is dedicated to author's}
\rightline{75-th anniversary\thanks{We present the article where the conception of so called "Kropina spaces" was introduced.
The article was published in Russian in "Trudy seminara po vektornomu i tenzornomu analizu" ("Workshops of the Seminar in vector and tensor Analysis"), vol. XI, 1961. 
Translated by P.D. Andreev. The interpreter tried to keep the style, notations and numerations of the author.}}

\vspace{2cm}

\author{V.K. Kropina}
\title[On projective two-dimensional Finsler spaces with special metric]{On projective two-dimensional Finsler spaces with special metric\\~\\ 
$L(x^\nu, X^v) = \frac{a_{\alpha\beta}(x^\nu)X^\alpha X^\beta}{b_\mu(x^\nu)X^\mu}$ and 
$L(x^\nu, X^v) = \sqrt[3]{a_{\alpha\beta\gamma}X^\alpha X^\beta X^\gamma}$}
\maketitle

Finsler space is called projective (or space with straight geodesics) if its geodesics are represented by linear equations in some special coordinate system (it is called projective). In his paper \cite{Be1} Berwald have solved the question of finding all projective two-dimensional Finsler spaces with metric

$$
L = \frac{(\alpha X + \beta Y)^2}{\gamma X + \delta Y}\eqno{({\rm I})}
$$

In this paper we solve the question of searching all projective two-dimensional Finsler spaces with metrics of more special type
$$
L=\frac{a_{11}(x,y) X^2 + 2 a_{12} XY + a_{22}Y^2}{b_1(x,y)X + b_2(x,y) Y}\eqno{(1)}\label{special}
$$
and indicate necessary and sufficient conditions for spaces with this metric to degenerate to Minkowski spaces. Also we establish that the curvature scalar of projective spaces with metric (1) cannot be constant (of course, with the exception of  Minkowski spaces)\footnote{Berwald in his paper \cite{Be1} showed that the same circumstance holds for projective spaces with metric (I)}.

The question of existence of two-dimensional projective Finsler spaces with metric given by the formula
$$
L^3 = A(x,y)X^3 + B(x,y) Y^3 + 3 C(x,y)X^2Y + 3D(x,y) XY^2\eqno{(1')}
$$
in assumption that the discriminant
$$
R = (AB-DC)^2 - 4(AD-C^2)(CB-D^2) \ne 0, \eqno{(1'')}
$$
is different from zero is testing in conclusion of the paper. We estimate that every space of such type is Minkowski space.

\textbf{1.} Coefficients $b_1$ and $b_2$ in expression of metric (1) can not be equal to zero simultaneously, hence we may assume without lost of generality that 
$$b_1 \ne 0\eqno{(2)}$$
and consequently the metric function (1) can be represented in following form
$$L = \frac{AX^2+BXY+CY^2}{X+DY}.\eqno{(3)}$$
We are assuming from now that second degree lines
$$L = \operatorname{const}\eqno{(4)}$$
are non-reducible, i.e.
$$\Delta = AD^2 - BD + C \ne 0. \eqno{(5)}$$
The necessary and sufficient conditions for Finsler space to have straight geodesics is the condition on its metric function to satisfy the following system of differential equations (cf. \cite{Be2}) in some special coordinate system (we call it projective):
\label{prcond}
$$\frac{\partial L}{\partial x^i} - \frac{\partial^2 L}{\partial x^\alpha\partial X^i}X^\alpha = 0. \eqno{(6)}$$

The system (6) in considered case (when metric function has form (3))\ can be reduced to the following system of differential equations for functions $A, B, C, D$:
$$
( {\rm I}) \left\{
			\begin{array}{lr}
				D\frac{\partial A}{\partial x^1} + \frac{\partial A}{\partial x^2} - \frac{\partial B}{\partial x^1} +A\frac{\partial D}{\partial x^1} = 0&\qquad (1_{\rm I})\\ & \\
				D^2\frac{\partial A}{\partial x^1} -AD\frac{\partial D}{\partial x^1} +3D \frac{\partial A}{\partial x^2} - 2\frac{\partial C}{\partial x^1} - D\frac{\partial B}{\partial x^1} +2B\frac{\partial D}{\partial x^1} = 0&\qquad (2_{\rm I})\\ & \\
				2D^2\frac{\partial A}{\partial x^2} - 3D \frac{\partial C}{\partial x^1} - \frac{\partial C}{\partial x^2} +D \frac{\partial B}{\partial x^2} +3C \frac{\partial D}{\partial x^1} + B \frac{\partial D}{\partial x^2} - 2AD \frac{\partial D}{\partial x^2} = 0 & \qquad (3_{\rm I})\\ & \\
				D^2\frac{\partial C}{\partial x^1} - D^2 \frac{\partial B}{\partial x^2} + D\frac{\partial C}{\partial x^2} -CD \frac{\partial D}{\partial x^1} + BD \frac{\partial D}{\partial x^1} + BD \frac{\partial D}{\partial x^2} -2C \frac{\partial D}{\partial x^2} = 0 & \qquad (4_{\rm I})
			\end{array}
	\right.
$$

\textbf{2.} Assuming at first that $D \ne 0$, we reorganize the system (I) as following. Subtracting equation $(1_{\rm I})$ multiplied by D from $(2_{\rm I})$, one get
$$D\frac{\partial A}{\partial x^2} - AD\frac{\partial D}{\partial x^1} - \frac{\partial C}{\partial x^1} +B\frac{\partial D}{\partial x^1} = 0. \eqno{(2'_{\rm I})}$$
Addition of equation $(4_{\rm I})$ with $(3_{\rm I})$ multiplied by $D$ gives
$$D^3\frac{\partial A}{\partial x^2} - D^2\frac{\partial C}{\partial x^1} + CD \frac{\partial D}{\partial x^1} + BD\frac{\partial D}{\partial x^2} - AD^2\frac{\partial D}{\partial x^2} - C \frac{\partial D}{\partial x^2} = 0. \eqno{(3'_{\rm I})}$$
Multiplying all items of $(2'_{\rm I})$ by $D^2$ and subtracting $(3_{\rm I})$ from the result, one finds
$$D\frac{\partial D}{\partial x^1}(AD^2-BD+C) - \frac{\partial D}{\partial x^2}(AD^2 - BD +C) = 0. \eqno{(3''_{\rm I})}$$

Dividing $(3''_{\rm I})$ by the expression $AD^2 - BD + C$ which is different from zero by (5), one gets
\[D \frac{\partial D}{\partial x^1} - \frac{\partial D}{\partial x^2} = 0. \eqno{(7)}\]
Hence system (I) can be changed by the following equivalent system (in frames of assumption $D \ne 0$)
\[
(\rm{I}')\left\{
			\begin{array}{lr}
				D\frac{\partial A}{\partial x^1} + \frac{\partial A}{\partial x^2} - \frac{\partial B}{\partial x^1} +A\frac{\partial D}{\partial x^1} = 0&\qquad (1_{\rm I'})\\ & \\
				D^2\frac{\partial A}{\partial x^1} -AD\frac{\partial D}{\partial x^1} +3D \frac{\partial A}{\partial x^2} - 2\frac{\partial C}{\partial x^1} - D\frac{\partial B}{\partial x^1} +2B\frac{\partial D}{\partial x^1} = 0&\qquad (2_{\rm I'})\\ & \\
				D \frac{\partial D}{\partial x^1} - \frac{\partial D}{\partial x^2} = 0 & \qquad (3_{\rm I'})\\ & \\
				D^2\frac{\partial C}{\partial x^1} - D^2 \frac{\partial B}{\partial x^2} + D\frac{\partial C}{\partial x^2} -CD \frac{\partial D}{\partial x^1} + BD \frac{\partial D}{\partial x^1} + BD \frac{\partial D}{\partial x^2} -2C \frac{\partial D}{\partial x^2} = 0 & \qquad (4_{\rm I'})
			\end{array}
	\right.
\]

\textbf{3.} First two equations of ($\rm{I}'$) can be expressed in following form (resolved relatively derivatives $\frac{\partial A}{\partial x^1}$ and $\frac{\partial A}{\partial x^2}$):
\[
(\alpha)\ \left\{
\begin{array}{lr}
				\frac{\partial A}{\partial x^1} = \frac1D \frac{\partial B}{\partial x^1} -  \frac{1}{D^2} \frac{\partial C}{\partial x^1} - \frac{2A}{D}\frac{\partial D}{\partial x^1} +\frac{B}{D^2}\frac{\partial D}{\partial x^1}&\hspace{6cm} (\alpha_1)\\ & \\
				\frac{\partial A}{\partial x^2} = \frac{1}{D}\frac{\partial C}{\partial x^1} + A \frac{\partial D}{\partial x^1} - \frac{B}{D}\frac{\partial D}{\partial x^1}.&\hspace{6cm} (\alpha_2)
			\end{array}
	\right.
\]

Let us compile conditions of integrability of the system $(\alpha)$. Differentiating the equation $(\alpha_1)$ by $x^2$ and taking into attention equations $(\alpha_2)$ and $(3_{\rm  I'})$, we get

\begin{align*}
\frac{\partial^2 A}{\partial x^1 \partial x^2} = - \frac1D \frac{\partial D}{\partial x^1}\frac{\partial B}{\partial x^1}& + \frac1D\frac{\partial^2 B}{\partial x^1\partial x^2} - \frac{1}{D^2}\frac{\partial^2 C}{\partial x^1\partial x^2} - \\
&-\frac{2A}{D}\frac{\partial^2 D}{\partial x^1\partial x^2} + \frac{1}{D^2}\frac{\partial B}{\partial x^2}\frac{\partial D}{\partial x^1} + \frac{B}{D^2}\frac{\partial^2 D}{\partial x^1\partial x^2}. \hspace{2.5cm} (9)
\end{align*}

After differentiating the equation $(3_{\rm I'})$ by $x^1$, we get
\[\frac{\partial^2 D}{\partial x^2 \partial x^1} = \left(\frac{\partial D}{\partial x^1}\right) + D \frac{\partial^2 D}{\partial {x^1}^2}.\eqno{(10)}\]

Differentiating the equation $(4_{\rm I'})$ by $x^1$, we find

\begin{align*}
\frac{\partial^2 C}{\partial x^2\partial x^1} = 2\frac{\partial C}{\partial x^1}\frac{\partial D}{\partial x^1} &+3C\frac{\partial^2 D}{\partial {x^1}^2} - D \frac{\partial B}{\partial x^1}\frac{\partial D}{\partial x^1} -\\
&-B\frac{\partial^2 D}{\partial x^2\partial x^1} + \frac{\partial D}{\partial x^1}\frac{\partial B}{\partial x^2} +D\frac{\partial^2 B}{\partial x^2\partial x^1} - D\frac{\partial^2 C}{\partial {x^1}^2}. \hspace{2.5cm} (11)
\end{align*}

Substituting the expression for $\frac{\partial^2 C}{\partial x^1\partial x^2}$ from (11) to equality (8) and taking into account (10), we get
\begin{align*}
\frac{\partial^2 A}{\partial x^2 \partial x^1} = - \frac{2}{D^2} \frac{\partial D}{\partial x^1}\frac{\partial C}{\partial x^1}& + \frac1D\frac{\partial^2 C}{\partial {x^1}^2} - \frac{1}{D^2}\frac{\partial^2 C}{\partial x^1\partial x^2} -\frac{2A}{D}\left(\frac{\partial D}{\partial x^1}\right)^2 + \\ 
&+ \frac{2B}{D^2}\left(\frac{\partial D}{\partial x^1}\right)^2 - \frac{B}{D}\frac{\partial^2 D}{\partial {x^1}^2} + A\frac{\partial^2 D}{\partial {x^1}^2}. \hspace{3.5cm} (12)
\end{align*}

As it follows from relations (12) and (9), the integrability condition
\[\frac{\partial^2 A}{\partial x^1\partial x^2} = \frac{\partial^2 A}{\partial x^2\partial x^1}\eqno(13)
\]
leads to the equality
\[3(AD^2 -BD +C)\frac{\partial^2 D}{\partial {x^1}^2} = 0, \eqno{(14)}
\]
i.e. in view of inequality (5), to equation
\[
\frac{\partial^2 D}{\partial {x^1}^2} = 0.
\]

So the function $D(x^1, x^2)$ is to satisfy following system of differential equations:
\[\hspace{5cm}
(\beta)\	\left\{
			\begin{array}{ll}
				\frac{\partial^2 D}{\partial {x^1}^2} & \hspace{6cm}(16)\\
				&\\
				\frac{\partial D}{\partial x^2} - D \frac{\partial D}{\partial x^1}& \hspace{6cm} (17)
 			\end{array}
		\right.
\]

The integrating of the system $(\beta)$ gives following solutions:

\[D = \frac{x^1 + k_1}{k_2 - x^2}, \eqno{(18)}\]
where $k_1$ and $k_2$ are arbitrary constants;
\[D = \operatorname{const} = k_1. \eqno{(19)}\]

\textbf{4.} Now we detect an essential geometric corollary of obtained result. The relation
\[X + DY = 0\eqno{(20)}\]
defines the direction (we will call it singular), tangent to curves $L = \operatorname{const}$  for every point\footnote{in the tangent centriaffine space; the contact is to be in the center of this space.}  (this is the unique direction for which the metric function is indefinite). It follows from founded expressions for D [(18) and (19)], that the flow lines for this field of directions (i. e. curves, tangents for which in every point have this singular direction) are characterized in our projective coordinate system by linear equations
\[x^1 - x^2k + k_1 + k_2k = 0 \eqno{(21)}\]
or
\[x^1 + k_1x^2 - k_1k = 0. \eqno{(22)}\]

It follows that these curves are "straights"\ (i.e. they maps onto straight lines in considering geodesic map of given space onto affine space). The "straightness"\ of the flow lines suggests the idea to pass to the new projective coordinate system where these lines will be coordinate lines $x_1 = \operatorname{const}$.

We show that such a choice of coordinate system is really possible. The singular direction has coordinates $\{D, -1\}$ in the initial coordinate system. The equation
\[\frac{\partial \varphi}{\partial x^1}D - \frac{\partial \varphi}{\partial x^2}\eqno{(25)}\]
is necessary and sufficient condition for the singular direction to have first coordinate equal to 0 in the new coordinate system $\bar x^1, \bar x^2$ which is connected with previous one by relations
\[\bar x^1 = \varphi(x^1, x^2) \eqno{(23)}\]
\[\bar x^2 = \psi(x^1, x^2). \eqno{(24)}\]
The integration of differential equation (25) gives:
\[\varphi = f\left(\frac{x^1 + k_1}{x^2 - k_2}\right) \eqno{(26)}\]
if $D$ is defined by relation (18) and gives also
\[\varphi = f(x^1 + k_1 x^2), \eqno{(27)}\]
if $D = \operatorname{const} = k_1$. Here $f$ is arbitrary function.

It follows namely from here that $D$ always may vanish in the result of projective (i.e. moving projective coordinate system to new projective one) coordinate transformation. It suffices to put $f(t)\equiv t$ and take
\[\psi = \frac{1}{x^2 - k_2} \eqno{(28)}\]
in the first case (when $\varphi$ is defined by relation (26)) or take
\[\psi = x^2 \eqno{(29)}\]
in the second case (when $\varphi$ is defined by relation (27)).

\textbf{5.} So we have shown that for every space with metric (3) with straight geodesics there exists a projective coordinate system for which $D = 0$.

Hence one may solve the system (1) putting $D=0$ for searching for all spaces with metric (3) and straight geodesics.

The system (I) becomes 
\[({\rm I}'')
		\left\{
			\begin{array}{l}
				\frac{\partial A}{\partial x^2} - \frac{\partial B}{\partial x^1} = 0\\
				\\
				\frac{\partial C}{\partial x^1} = 0\\
				\\
				\frac{\partial C}{\partial x^2} = 0,
			\end{array}
		\right.
\]
and it follows
\[
\hspace{5cm}
		\left\{
			\begin{array}{ll}
				C=k = \operatorname{const}& \hspace{6cm} (30)\\
				\\
				A = \frac{\partial \varphi}{\partial x^1}&\hspace{6cm} (31)\\
				\\
				B = \frac{\partial \varphi}{\partial x^2},&\hspace{6cm} (32)
			\end{array}
		\right.
\]
where $\varphi(x^1, x^2)$ is an arbitrary function of $x^1$, $x^2$.

Since we have
\[\Delta = k\eqno{(33)}\]
in this case,  the necessary and sufficient condition for holding the inequality (5) is
\[k \ne 0. \eqno{(34)}\]
So we have achieved following expression (in the special projective coordinate system) for the metric function (3) of the space with straight geodesics:
\[L = \frac{\frac{\partial \varphi}{\partial x^1} X^2 + \frac{\partial \varphi}{\partial x^2}XY + kY^2}{X}. \eqno{(35)}\]
After coordinate transformation
\[\bar x^1 = k x^1, \eqno{(36)}\]
\[\bar x^2 = k x^2, \eqno{(37)}\]
we get finally
\[L = \frac{\frac{\partial \varphi}{\partial x^1} X^2 + \frac{\partial \varphi}{\partial x^2}XY + Y^2}{X}. \eqno{(38)}\]
So we showed that Finsler space with metric (3) has straight geodesics iff its metric function in some coordinate system has view (38), where $\varphi(x^1, x^2)$ is an arbitrary function.

R e m a r k. In view of correlation of our result with the result of Berwald, we detect spaces with metrics of parabolic type (i. e. such spaces that lines $L = \operatorname{const}$ are parabolas) among founded spaces.

The necessary and sufficient condition for space with metric function (38) to satisfy such property is the condition for the function $\varphi$ to be solution of the differential equation
\[4 \frac{\partial \varphi}{\partial x^1} - \left(\frac{\partial \varphi}{\partial x^2}\right)^2 = 0. \eqno{(39)}\]
It follows from here that the function $\varphi$ must have view
\[\varphi = \varkappa_1x^2 + \frac{\varkappa_1^2}{4} x^1 + \varkappa_2, \eqno{(40)}\]
where $\varkappa_1$, $\varkappa_2$ are constants, or
\[\varphi = x^2 a + \frac{a^2}{4}x^1 + \omega(a), \eqno{(41)}\]
where $\omega$ --- arbitrary function of one parameter, and $a(x^1, x^2)$ --- function, defined from equation
\[x^2 + \frac{a}{2}x^1 + \omega'(a) = 0.\eqno{(42)}\]
We note that the expression (40) is out of interest because it leads to Minkowski space.
From equality (41) we find:
\[\frac{\partial \varphi}{\partial x^2} = a + x^2 \frac{\partial a}{\partial x^2} + \frac{a}{2} \frac{\partial a}{\partial x^2}x^1 + \omega'(a)\frac{\partial a}{\partial x^2}. \eqno{(43)}\]
Taking into account relation (42) we find from this:
\[\frac{\partial \varphi}{\partial x^2} = a(x^1, x^2). \eqno{(44)}\]
So we obtain the result: the metric function (3) of parabolic type in the case of the space with straight geodesics has in some coordinate system following view:
\[L = \frac{(a X + Y)^2}{X}, \eqno{(45)}\]
where $a(x^1, x^2)$ --- function defined from the equation
\[x^2 + x^1a + \sigma(a) = 0. \eqno{(46)}\]
Here $\sigma$ is an arbitrary function of one variable.

So we have recovered the result obtained by Berwald [cf. \cite{Be1}, formulas (16.18) and (16.20)].

\textbf{6.} Naturally arising question is in what case the space with metric function (38) is Minkowskian? What kind must be the function $\varphi(x^1, x^2)$ in this case?

One can show that the necessary and sufficient condition for the space with straight geodesics to be Minkowskian consists in holding (in projective coordinate system) following system of differential equations (we omit the proof here):
\[({\rm II})\
	\left\{
		\begin{array}{l}
			p_{ij} = 0\\
			\\
			\frac{\partial p_j}{\partial x^i} - p_jp_i = 0,
		\end{array}
	\right.
\]
where
\[p(x^\alpha, X^\alpha) = \frac{\frac{\partial L}{\partial x^\nu}X^\nu}{2L}, \eqno{(47)}\]
\[p_i = \frac{\partial p}{\partial X^i}, \quad p_{ij} = \frac{\partial^2 p}{\partial X^i\partial X^j}. \eqno{(48)}\]

After computations (which are omitted because of their cumbersomity), we may verify that the system (II) in considering case can be reduced to following system of differential equations for the function $\varphi(x^1, x^2)$:
\[
(\rm{II}')\ 
	\left\{
		\begin{array}{l}
			2 \frac{\partial^2 \varphi}{\partial x^1\partial x^2} - \frac{\partial \varphi}{\partial x^1}\frac{\partial^2 \varphi}{(\partial x^1)^2} = 0\\
			\\
			\frac{\partial^2 \varphi}{(\partial x^2)^2} - \frac{\partial \varphi}{\partial x^2}\frac{\partial^2\varphi}{(\partial x^1)^2} = 0\\
			\\
			\frac{\partial^3\varphi}{(\partial x^1)^3} = 0.
		\end{array}
	\right.
\]
In the result of integration of the system ($\rm{II}'$) we find that 
\[\varphi = k_1 x^1 + k_2 x^2 + k_2 \eqno{(49)}\]
or
\[\varphi = \frac{1}{k_1 - x^2}[(x^1)^2 + k_2x^1 + k_3] + k_4. \eqno{(50)}\]
So we conclude: the necessary and sufficient condition for the space with metric function (38) to be Minkowskian is the condition on the function $\varphi$ in equation (38) to be of type defined by equality (49) or (50).

\textbf{7.} The result obtained above lets us to clear up, if there exist spaces which are not Minkowski among founded spaces (38) of constant curvature with straight geodesics.

It is necessary for Finsler space to be a space of constant curvature, that its curvature tensor $K_{ijkl}$ is skew-symmetric in both pairs of indices $i,j$ and $k,l$. This is a corollary of formulas detected by Berwald [cf. \cite{Be3}, formula (64) and \cite{Be4}, formula (13.7)]

Now we show that the relation (in projective coordinate system)
\[\frac{\partial p_i}{\partial x^j} - \frac{\partial p_j}{\partial x^i} = 0 \eqno{(51)}\]
follows from the skew symmetry (in both pairs of indices) of the curvature tensor of Finsler space with straight geodesics.

As it is known \cite{Be1}, the connection $G^i_{jk}$ of Finsler space with straight geodesics in the projective coordinate system has following view:
\[G^i_{jk} = \delta^i_jp_k + \delta^i_kp_j + X^ip_{jk}. \eqno{(52)}\]
From this we obtain after bulky calculations following expression for the curvature tensor $K^i_{jkl}$:
\[K^i_{jkl} = \delta^i_j(P_{kl} - P_{lk}) + \delta^i_k P_{jl} + X^i\left(\frac{\partial P_{kl}}{\partial X^j} - \frac{\partial P_{lk}}{\partial X^j}\right), \eqno{(53)}\]
where
\[P_{ij} = \frac{\partial p_i}{\partial x^j} - p_ip_j - p_{ij}p. \eqno{(54)}\]
If the curvature tensor $K_{ijkl}$ is symmetric by indices $i$ and $j$, then
\[K^\alpha_{\alpha kl} = 0. \eqno{(55)}\]
But it follows from the equality (53), taking into account the 0-homogeneouity of the function $P_{ij}$ in $X^\alpha$, that
\[K^\alpha_{\alpha kl} = (n+1)(P_{kl} - P_{lk}). \eqno{(56)}\]
Comparing with (55), we have
\[P_{kl} - P_{lk} = 0. \eqno{(57)}\]
So, holding (in the projective coordinate system) the system of differential equations (58)  is necessary for the Finsler space with straight geodesics to be a space of constant curvature.

\textbf{8.} Returning to the question asked in the beginning of p. \textbf{7}, we rewrite the system (58) in considering case, when the metric function has view (38); we get the following system of differential equations for the function $\varphi$:

\[\frac{\partial^3 \varphi}{(\partial x^1)^3} = 0, \eqno{(1^1)}\]
\[2 \frac{\partial^3 \varphi}{(\partial x^1)^2\partial x^2} - \left[\frac{\partial^2 \varphi}{(\partial x^1)^2}\right]^2 = 0, \eqno{(2^1)}\]
\[\frac{\partial\varphi}{\partial x^2}\frac{\partial^3 \varphi}{(\partial x^1)^2\partial x^2} - 2 \frac{\partial^2 \varphi}{(\partial x^1)^2}\frac{\partial^2 \varphi}{\partial x^1\partial x^2} + \frac{\partial^3 \varphi}{(\partial x^2)^2\partial x^1} = 0, \eqno{(3^1)}\]
\begin{align*}
\hspace{2cm}
3\frac{\partial \varphi}{\partial x^2}\left(\frac{\partial^2 \varphi}{\partial {x^1}^2}\right)^2 - 6\left(\frac{\partial \varphi}{\partial x^1}\right)^2\frac{\partial^3 \varphi}{\partial {x^1}^2\partial x^2} &+ 12 \frac{\partial^2\varphi}{\partial {x^1}^2}\frac{\partial^2 \varphi}{\partial x^1\partial x^2}\frac{\partial \varphi}{\partial x^1} -\\
- 12\left(\frac{\partial^2 \varphi}{\partial x^1\partial x^2}\right)^2 - 6 \frac{\partial^2 \varphi}{\partial {x^1}^2}\frac{\partial^2 \varphi}{\partial {x^2}^2} &+ 2\frac{\partial^3 \varphi}{\partial {x^1}^3}, \hspace{4cm} (4^1)
\end{align*}
\begin{align*}
\hspace{2cm}
5\frac{\partial \varphi}{\partial x^1}\frac{\partial \varphi}{\partial x^2} \frac{\partial^3 \varphi}{\partial {x^1}^2\partial x^2} - &14 \frac{\partial \varphi}{\partial x^2}\frac{\partial^2 \varphi}{\partial {x^1}^2}\frac{\partial^2 \varphi}{\partial x^1\partial x^2} + 14 \frac{\partial^2 \varphi}{\partial {x^2}^2}\frac{\partial^2 \varphi}{\partial x^1\partial x^2} - \\
- 2 \left(\frac{\partial \varphi}{\partial x^1}\right)^2\frac{\partial^3 \varphi}{\partial {x^2}^2\partial x^1} + &4 \frac{\partial \varphi}{\partial x^1}\left(\frac{\partial^2 \varphi}{\partial x^1\partial x^2}\right)^3 +2 \frac{\partial \varphi}{\partial x^1} \frac{\partial^2 \varphi}{\partial {x^1}^2}\frac{\partial^2 \varphi}{\partial {x^2}^2} - \\
- &3 \frac{\partial \varphi}{\partial x^1}\frac{\partial^3 \varphi}{\partial {x^2}^3} = 0,
\hspace{6cm} (5^1)
\end{align*}
\begin{align*}
\hspace{2cm}
2\left(\frac{\partial \varphi}{\partial x^2}\right)^2\frac{\partial^3 \varphi}{\partial {x^1}^2\partial x^2} - &2\frac{\partial \varphi}{\partial x^2}\frac{\partial^3 \varphi}{\partial {x^2}^3} - \left(\frac{\partial \varphi}{\partial x^1}\right)^2\frac{\partial^3 \varphi}{\partial {x^2}^3} -\\
- 2\frac{\partial \varphi}{\partial x^2}\frac{\partial^2 \varphi}{\partial {x^1}^2}\frac{\partial^2 \varphi}{\partial {x^2}^2} - 4 &\frac{\partial \varphi}{\partial x^2}\left(\frac{\partial^2 \varphi}{\partial x^1\partial x^2}\right)^2 + 6 \frac{\partial \varphi}{\partial x^1}\frac{\partial^2 \varphi}{\partial {x^2}^2}\frac{\partial^2 \varphi}{\partial x^1\partial x^2} +\\
+4\left(\frac{\partial^2 \varphi}{\partial {x^1}^2}\right)^2& - \frac{\partial \varphi}{\partial x^1}\frac{\partial \varphi}{\partial x^2}\frac{\partial^3 \varphi}{\partial {x^2}^2\partial x1},\hspace{4.5cm} (6^1)
\end{align*}
\begin{align*}
\hspace{2cm}
2\frac{\partial \varphi}{\partial x^2}\frac{\partial^2 \varphi}{\partial {x^2}^2}\frac{\partial^2 \varphi}{\partial x^1\partial x^2} &- \left(\frac{\partial \varphi}{\partial x^2}\right)^2\frac{\partial^3 \varphi}{\partial {x^2}^2\partial x^1} + \frac{\partial \varphi}{\partial x^1}\frac{\partial \varphi}{\partial x^2}\frac{\partial^3 \varphi}{\partial {x^2}^3} -\\
&-2\frac{\partial \varphi}{\partial x^1}\left(\frac{\partial^2 \varphi}{\partial {x^2}^2}\right)^2 = 0. \hspace{4.5cm} (7^1)
\end{align*}

It follows from the equation $(1^1)$ that the function $\varphi$ is of type
\[\varphi(x^1, x^2) = \alpha(x^2) (x^1)^2 + \beta(x^2)x^1 + \gamma(x^2). \eqno{(59)}\]
Substituting expressions for derivatives from this to equation $(2^1)$, we obtain differential equation for the function $\alpha(x^2)$:
\[\frac{\partial \alpha}{\partial x^2} - \alpha^2 =0. \eqno{(60)}\]
So we have
\[\alpha = \frac{1}{k_1 - x^2} \eqno{(61)}\]
or
\[\alpha = 0. \eqno{(62)}\]
Substititing  expressions for derivatives of the function $\varphi$ from (59) into the equation $(3^1)$ and taking into account expression (61), we find differential equation for the function $\beta(x^2)$:
\[\frac{\partial^2 \beta}{\partial {x^2}^2}(k_1 - x^2)^2 - 4(k_1 - x^2)\frac{\partial \beta}{\partial x^2} + 2\beta = 0. \eqno{(63)}\]
After integration we have
\[\beta = \frac{k_2x^2 + k_3}{(k_1 - x^2)^2}. \eqno{(64)}\]
When $\alpha = 0$ the equation $(3^1)$ gives us
\[\beta = \varkappa_1x^2 + \varkappa_2, \eqno{(65)}\]
where $\varkappa_1$ and $\varkappa_2$ are arbitrary constants.

Substititing  expressions for derivatives of the function $\varphi$ from (59) into the equation $(4^1)$ and taking into account equalities (61) and (64), we find differential equation for the function $\gamma(x^2)$:

\begin{align*}
\hspace{2cm}
(k_1 - x^2)\frac{\partial^3 \gamma}{\partial {x^2}^3} - 6(k_1 - &x^2)^5\frac{\partial^2 \gamma}{\partial {x^2}^2} + 6 (k_1 - x^2)^4 \frac{\partial \gamma}{\partial x^2} -
\\ -& 6(k_1k_2 + k_3)^2 = 0. \hspace{5cm} (66)
\end{align*}
Resolving this equation, we find:
\[\gamma=\frac{k_4}{2(k_1 - x^2)^2)} + \frac{k_5}{k_1 - x^2} + \frac{(k_1k_2+k_3)^2}{(k_1-x^2)^3} + k_6. \eqno{(67)}\]

If $\alpha$ and $\beta$ are defined by equalities (62) and (65) correspondingly, then equation $(4^1)$ leads to the expression
\[\gamma = (\varkappa_1)^2(x^2)^3 + \frac{\varkappa_3}{2} (x^2)^2 + \varkappa_4x^2 + \varkappa_5. \eqno{(68)}\]
The equation $(5^1)$ after substituting of derivatives founded from equality (59) taking into account (61), (64) and (67), leads to relation
\[k_1k_2 + k_3 = 0. \eqno{(69)}\]

In the case when $\alpha$, $\beta$ and $\gamma$ are defined by equalities (62), (65) and (68) correspondingly, the equation $(5^1)$ gives
\[\varkappa_1 = 0.\eqno{(70)}\]

Manipulating with the equation $(5^1)$ similarly as with previous and taking into account relation (69), we find
\[k_4 = 0, \eqno{(71)}\]
when $\alpha$, $\beta$ and $\gamma$ are defined by relations (61), (64) and (67) and
\[\varkappa_3 = 0, \eqno{(72)}\]
if  $\alpha$, $\beta$ and $\gamma$ are defined by (62), (65) and (68) correspondingly and the equality (70) holds.

In frames of relations (61), (64), (67), (69) and (71) equality (59) becomes of form
\[\varphi(x^1, x^2) = \frac{1}{k_1 - x^2}[ (x^1)^2 - k_2x^1 + k_5] + k_6, \eqno{(73)}\]
and it follows from equalities (62), (65), (68), (70) and (72):
\[\varphi = \varkappa_2 x^1 + \varkappa_4 x^2 + \varkappa_5.\eqno{(74)}\]
But as it was mentioned above (p. \textbf{7}), the space with metric (38) and function $\varphi$ of type (73) or (74) is Minkowskian. So, we proved that spaces with metric function (38) can be of constant curvature if and only if they are Minkowski spaces.

\textbf{9.} Next we will turn to the question on projectivity of Finsler spaces with the metric given by formula $(1')$. Using the projectivity condition (6) from p. \pageref{prcond}
\[\frac{\partial L}{\partial x^i} - \frac{\partial^2 L}{\partial x^\alpha \partial X^i}X^\alpha\]
and after necessary calculations\footnote{We omit calculations; we note only that we use inequality $(1'')$ in the calculation process.} we obtain following system of differential equations for functions $A(x^1, x^2)$, $B(x^1, x^2)$, $C(x^1, x^2)$, $D(x^1, x^2)$ of Finsler space $(1')$ with straight geodesics:
\[
(\alpha)\ \left\{
			\begin{array}{l}
				\frac{\partial A}{\partial x^1} = \frac{2A}{B} \frac{\partial B}{\partial x^1} - \frac{AD}{B^2}\frac{\partial B}{\partial x^2}\\
				\\
				\frac{\partial A}{\partial x^2} = \frac{AB-CD}{2B^2}\frac{\partial B}{\partial x^2} + \frac{C}{B}\frac{\partial B}{\partial x^1}\\
				\\
				\frac{\partial C}{\partial x^1} = \frac{5C}{3B}\frac{\partial B}{\partial x^1} + \frac{AB - 5CD}{\frac{\partial B}{\partial x^2}}\\
				\\
				\frac{\partial C}{\partial x^2} = \frac{2D}{3B}\frac{\partial B}{\partial x^1} + \frac{2BC- D^2}{3B^2}\frac{\partial B}{\partial x^2}\\
				\\
				\frac{\partial D}{\partial x^1} = \frac{4D}{3B}\frac{\partial B}{\partial x^1} + \frac{BC - 2 D^2}{3B^2}\frac{\partial B}{\partial x^2}\\
				\\
				\frac{\partial D}{\partial x^2} = \frac13 \frac{\partial B}{\partial x^1} + \frac23 \frac{D}{B}\frac{\partial B}{\partial x^2}\\
				\\
				6B \frac{\partial^2 B}{\partial {x^2}^2} = 7\left(\frac{\partial B}{\partial x^2}\right)^2\\
				\\
				\frac{\partial^2 B}{\partial x^1\partial x^2} = \frac{4}{3B}\frac{\partial B}{\partial x^1}\frac{\partial B}{\partial x^2} - \frac16 \frac{D}{B^2}\left(\frac{\partial B}{\partial x^2}\right)^2\\
				\\
				\frac{\partial^2 B}{\partial {x^1}^2} = \frac{4}{3B} \left(\frac{\partial B}{\partial x^1}\right)^2 + \frac{BC-2D^2}{6B^2}\left(\frac{\partial B}{\partial x^2}\right)^2
			\end{array}
		\right.
\]
and all calculations are carried in assumption that $B \ne 0$. The case $B = 0$ is left to special consideration.

\textbf{10.} It is worthwhile in this stage, i.e. before integration of the system $(\alpha)$, to clear up, whether there exists among projective Finsler spaces with metric (1), spaces which are essentially Finsler, that is spaces which are not Minkowski. For this we use conditions (II) found out in p. \textbf{6}:
\[p_{ij} = 0, \frac{\partial p_j}{\partial x^i} - p_ip_j = 0.\]
Taking into attention first degree homogeneouity of the function $p(x^\alpha, X^\alpha)$ in variables $X^\alpha$, we easily convince ourselves that the system (II) is equivalent to the following system
\[(\rm{III})\ 
	\left\{
		\begin{array}{l}
			p_{ij} = 0\\
			\\
			\frac{\partial p}{\partial x^i} = pp_i.
		\end{array}
	\right.
\]

Using formula (47) we find the function $p$ for considering Finsler space $(1')$ and after necessary calculations using the system of differential equations $(\alpha)$, we prove that function $p$ satisfies to the system of differential equations (III); consequently all spaces with metric $(1')$ with straight geodesics and satisfying to condition $B \ne 0$ (in special coordinate system, which is discussed here) are necessary Minkowski spaces.

\textbf{11.} It is left to consider the exceptional case.

Putting $B = 0$ (and $D \ne 0$)\footnote{The case $B=0$, $D=0$ is excluded by the inequality $(1'')$.}, after integrating of the system (6) we obtain the following expressions for functions $A$, $B$, $C$, $D$ in Finsler space with straight geodesics:
\[A = \frac{k_1}{(k_1x^1 + k_2)^6}(3k_1(x^2)^2 - 3 k_3 x^2 + k_4),\]
\[B = 0,\]
\[C = \frac{1}{(k_1x^1 + k_2)^5}(-2 k_1x^2 + k_3),\]
\[D = \frac{1}{(k_1x^1 + k_2)^4},\]
where $k_1$, $k_2$, $k_3$, $k_4$ are arbitrary constants. After straight calculations we verify that function $p$ satisfies to the system (III) and consequently in considering case  all Finsler spaces ($1'$) with straight geodesics and $B = 0$ (in special coordinate system) are Minkowski spaces.

Finally, we have cleared up that Finsler spaces with metric function satisfying the requirement $(1'')$ can have straight geodesics when and only when they are Minkowski spaces.

\begin{flushright}
\textit{The original paper is dated by 15.05.1958.}
\end{flushright}  


\begin{thebibliography}{99}

\bibitem{Be1} B e r w a l d\ \ L. 
On Finsler and Cartan Geometrie. III
Ann. of Math, 1941, {\bf 42}, 1, 84--112.


\bibitem{Be2} B e r w a l d\ \ L.
\"Uber die $n$-dimensionalen Geometrien konstanter Kr\"ummung, in denen die Geraden die K\"urzesten sind. Math. Z., 1929, \textbf{30}, \no 3, 449--469.

\bibitem{Be3} B e r w a l d\ \ L.
Untersuchung der Kr\"ummung allgemeiner metrischer R\"aume auf Grund des ihnen herrschenden Parallelismus. Math. Z., 1926, \textbf{25}, 1, 40--73

\bibitem{Be4} B e r w a l d\ \ L.
\"Uber Finslersche und Cartansche Geometrie. IV --- Ann. of Math., 1947, \textbf{48}, 3, 755--781.
 
 
\end{thebibliography}
\end{document}